\pdfoutput=1
\RequirePackage{ifpdf}
\ifpdf 
\documentclass[pdftex]{sigma}
\else
\documentclass{sigma}
\fi

\numberwithin{equation}{section}
\newtheorem{Theorem}{Theorem}[section]
\newtheorem{Corollary}[Theorem]{Corollary}
{ \theoremstyle{definition}
\newtheorem{Note}[Theorem]{Note}
\newtheorem{Remark}[Theorem]{Remark} }

\def\CC{\mathbb{C}}
\def\RR{\mathbb{R}}
\def\NN{\mathbb{N}}
\def\ZZ{\mathbb{Z}}
\def\bfy{\mathbf{y}}
\def\bfs{\mathbf{s}}

\begin{document}
\allowdisplaybreaks

\newcommand{\arXivNumber}{1804.11273}

\renewcommand{\thefootnote}{}

\renewcommand{\PaperNumber}{117}

\FirstPageHeading

\ShortArticleName{Truncated Solutions of Painlev\'e Equation ${\rm P}_{\rm V}$}

\ArticleName{Truncated Solutions of Painlev\'e Equation $\boldsymbol{{\rm P}_{\rm V}}$\footnote{This paper is a~contribution to the Special Issue on Painlev\'e Equations and Applications in Memory of Andrei Kapaev. The full collection is available at \href{https://www.emis.de/journals/SIGMA/Kapaev.html}{https://www.emis.de/journals/SIGMA/Kapaev.html}}}

\Author{Rodica D.~COSTIN}

\AuthorNameForHeading{R.D.~Costin}

\Address{The Ohio State University, 231 W 18th Ave, Columbus, OH 43210, USA}
\Email{\href{mailto:costin.10@osu.edu}{costin.10@osu.edu}}
\URLaddress{\url{https://people.math.osu.edu/costin.10/}}

\ArticleDates{Received May 01, 2018, in final form October 25, 2018; Published online October 31, 2018}

\Abstract{We obtain convergent representations (as Borel summed transseries) for the five one-parameter families of truncated solutions of the fifth Painlev\'e equation with nonzero parameters, valid in half planes, for large independent variable. We also find the position of the first array of poles, bordering the region of analyticity. For a special value of this parameter they represent tri-truncated solutions, analytic in almost the full complex plane, for large independent variable. A brief historical note, and references on truncated solutions of the other Painlev\'e equations are also included.}


\Classification{33E17; 34M30; 34M25}

\renewcommand{\thefootnote}{\arabic{footnote}}
\setcounter{footnote}{0}
\vspace{-2mm}

\section{Introduction}

\subsection{Historical notes} By the middle of the 19th century it became apparent that solutions of many linear differential equations should be considered new, ``special'' functions. The natural question then arose: can {\em nonlinear} differential equations have solutions that could be thought of being special functions? Fuch's intuition was that the answer is affirmative if ``solutions have only fixed branch points none of which depend on the initial conditions''~\cite{story}. This is now called {\em the Painlev\'e property}, the fact that solutions of an equation are meromorphic on a common Riemann surface. Fuchs studied first order equations having this property, concluding in 1884 that all such equations can be solved in terms of previously known functions; his results were later extended by Poincar\'e~\cite{story}.

The idea that absence of movable branch points would mean integrability was then used by Sofie Kowalevski in the study of the rotation of a solid about a fixed point, for which she disco\-ve\-red a third integrable case~\cite{Sofia} (previous two cases being disco\-ve\-red by Euler and Lagrange), a~discovery for which she was awarded the Prix Bordin of the French Academy of Science in 1888.
Around the turn of the 20th century, Painlev\'e, Picard and Gambier studied nonlinear second order differential equations, rational in~$y$ and~$y'$, discovering that those possessing the Painlev\'e property can be brought to fifty canonical forms. Of these fifty equations, all but six could be solved in terms of earlier known functions. Painlev\'e went on to show that generic solutions of the remaining six equations cannot be expressed in terms of earlier known functions, or in terms of each other~\cite{Painleve}. These equations are now known as the Painlev\'e equations, denoted ${\rm P}_{\rm I}$ up to~${\rm P}_{\rm VI}$, and their solutions as the {\em Painlev\'e transcendents}.

Having been discovered as a result of a purely theoretical inquiry, the Painlev\'e transcendents have appeared later in modern geometry, integrable systems~\cite{Ablowitz}, statistical mecha\-nics~\cite{7, 8,6}, and recently in quantum field theory. Their practical importance makes it necessary to study their properties in detail and to develop good numerical methods for their calculation~\cite{Clarkson}.

\subsection{Truncated solutions of Painlev\'e equations}
\looseness=-1 It is known that generic Painlev\'e transcendents have poles in any sector towards infinity. But there are special solutions, called tronqu\'ee (or truncated) which are free of poles in some sectors, at least for large values of the independent variable. It was later conjectured that truncated solutions have no poles whatsoever in these sectors, the Dubrovin--Novokshenov conjecture~\cite{Dubrovin,Novo}.

For ${\rm P}_{\rm I}$, Boutroux showed that there are five special sectors in the complex plane, each of opening $2\pi/5$, where solutions may lack poles. There are truncated solutions free of poles in two adjacent such sectors, for large~$x$. Among these, there are solutions free of poles in four sectors~-- tri-truncated solutions~\cite{Boutroux}.

In each sector, there is a one parameter family of truncated solutions, all asymptotic to the same power series; these solutions can be distinguished by a parameter multiplying a small exponential; a complete formal solution can be obtained as a transseries, mixing negative powers of the variable and exponentially small terms. The form of the exponentially small terms as well as the location of the first array of poles of truncated solutions were found in \cite{P1us}. The Dubrovin--Novokshenov conjecture was proved for ${\rm P}_{\rm I}$ \cite{Dconj}.

The Stokes constant (which can be viewed as a complex number, attached to the equation, controlling the relation between the exponentially small terms associated to different sectors) was calculated for the tritronqu\'ee solution for ${\rm P}_{\rm I}$ for the first time by Kapaev~\cite{Kapaev1}, see also~\cite{Kapaev}, using the isomonodromy method. Kapaev then obtained the complete description of the global asymptotic behavior of the tronqu\'ee solutions of~${\rm P}_{\rm I}$ and connection formulae~\cite{Kapaev2}. The Stokes constant was re-calculated later using WKB methods~\cite{Takei}, and again later by continuation of a~tri-truncated solution through the sector with poles~\cite{OCRDC}.

More recently, Kapaev considered the study of the tronqu\'ee solution of the ${\rm P}_{\rm I2}$ equation (the second member of the ${\rm P}_{\rm I}$ hierarchy); a detailed global asymptotic analysis of the trintronqu\'es solutions is found in~\cite{Kapaev6}.

For ${\rm P}_{\rm II}$ Boutroux showed that there are six special sectors, with truncated solutions free of poles in two adjacent sectors, and tri-truncated solutions free of poles in four sectors.
Using the Riemann--Hilbert approach, Kapaev gave a complete description of the global asymptotic behavior of the tronqu\'ee solutions of ${\rm P}_{\rm II}$, together with all relevant connection formulae~\cite{Kapaev4,Kapaev3}, see also Chapter~11 of the monograph~\cite{Fokas}.
Existence of tri-truncated solutions of the ${\rm P}_{\rm II}$ hierarchy was shown in~\cite{Nalini}. The Dubrovin--Novokshenov conjecture was established for the Hastings--McLeod solution (tri-truncated in pairs of non-adjacent sectors)~\cite{Min}.

For~${\rm P}_{\rm IV}$, a~quite detailed global asymptotic analysis of its tronqu\'ee solutions was obtained in Kapaev's work \cite{Kapaev5}. The connection formulae were found by Its and Kapaev~\cite{IK}, based on the Riemann--Hilbert isomonodromy method.
Existence of truncated solution fo~${\rm P}_{\rm III}$ and~${\rm P}_{\rm IV}$ was re-established in \cite{Lin} following methods in~\cite{Kitaev}, and using a different method in~\cite{Xia}, where the location of the first array of poles was also found.
An overview of~${\rm P}_{\rm VI}$ is contained in \cite{Guzzetti}; see also \cite{WePvi}.
The truncated solutions of the fifth Painlev\'e equation are the subject of the present article. It is to be noted that fixed singularities of ${\rm P}_{\rm V}$ can only be located at $0$ and $\infty$.

\subsection[Truncated solutions of Painlev\'e equation ${\rm P}_{\rm V}$]{Truncated solutions of Painlev\'e equation $\boldsymbol{{\rm P}_{\rm V}}$}
The fifth Painlev\'e equation{\samepage
\begin{gather}
{\rm P}_{\rm V}(\alpha,\beta,\gamma,\delta)\colon\nonumber \\
w''= \left( \frac 1{2w}+ \frac 1{ w -1} \right) {w' }^{2}-{\frac {w' }{x}}+{\frac { (w-1)
^{2}}{{x}^{2}} \left( \alpha w +{\frac {\beta}{w }} \right) }+{\frac {\gamma w }{x}}+{\frac {\delta w (w +1) }{w -1}} \label{P5}
\end{gather}
is known to be reducible to ${\rm P}_{\rm III}$ if $\delta=0$ \cite{Gromak}.}

We assume $\alpha\beta\delta\ne 0$. The following algebraic behaviors towards infinity are then possible for solutions of \eqref{P5} \cite{Parusnikova}:
\begin{alignat*}{4} 
& {\rm I.}\quad & &w=\pm \sqrt{\frac{\beta}{\delta}} x^{-1}+O\big(x^{-2}\big), \qquad&& x\to\infty, & \\
&{\rm II.} \quad && w=\pm\sqrt{\frac{\delta}{-\alpha}} x+O(1), \qquad &&x\to\infty, & \\
&{\rm III.} \quad && w= -1+O\big(x^{-1}\big), \qquad && x\to\infty.&
 \end{alignat*}

It is known that these five families represent asymptotic behaviors of truncated solutions, analytic in (almost) a half plane for large $|x|$; the position of the half plane is determined by the exponentially small terms \cite{11,Shimomura}.

In the present paper we express the exponentially small terms using a full formal solution (transseries), Borel summable to actual solutions. This yields one-parameter families of truncated solutions as series which converge in appropriate (almost) half-planes and large~$|x|$. For special values of the parameter we obtain tri-truncated transcendents.
Moreover, we find the location of the first array of poles beyond the sector of analyticity.

We use techniques used before in~\cite{P1us}. For completeness, the statements of the theorems used here are included in~Appendix~\ref{AppB}.

\subsection{Relations between different truncated solutions}

\begin{Remark}[\cite{Gromak}]\label{Rem1}\quad\begin{enumerate}\itemsep=0pt
\item[(i)] If $w(x)$ satisfies ${\rm P}_{\rm V}(\alpha,\beta,\gamma,\delta)$, then $1/w(x)$ satisfies ${\rm P}_{\rm V}(-\beta,-\alpha,-\gamma,\delta)$.
\item[(ii)] If $w(x)$ solves ${\rm P}_{\rm V}(\alpha,\beta,\gamma,\delta)$ then $w(x/\lambda )$ solves ${\rm P}_{\rm V}\big(\alpha,\beta,\gamma\lambda,\delta\lambda^2\big)$, for any $\lambda\ne 0$.
\end{enumerate}
 \end{Remark}

By Remark~\ref{Rem1}(i), truncated solutions in the family~II. are obtained as reciprocals of truncated of the family~I.

By Remark~\ref{Rem1}(ii), if $w(x)$ solves ${\rm P}_{\rm V}(\alpha,\beta,\gamma,\delta)$ with $w(x)\sim \sqrt{\frac{\beta}{\delta}} x^{-1}$, $x\to\infty$, then $w(-x)$ solves ${\rm P}_{\rm V}(\alpha,\beta,-\gamma,\delta)$ and satisfies $w(-x)\sim- \sqrt{\frac{\beta}{\delta}} x^{-1}$, $x\to\infty$.

We can assume any nonzero value for $\delta$, by rescaling $x$ and using Remark~\ref{Rem1}(ii).

It is interesting to note that there are other B\"acklund transformations as explained in \cite[Theorem~39.2]{Gromak}. These transform truncated solutions of the first two families into truncated solutions of the first two families, and truncated solutions in the family~III into solutions in the same family.

Due to these relations, it suffices to obtain results for the truncated Painlev\'e transcendents satisfying
\begin{gather*}
{\rm I_0} \quad w= \sqrt{{\beta}/{\delta}} x^{-1}+O\big(x^{-2}\big), \qquad x\to\infty, \qquad \text{for} \quad \delta=-1/2
\end{gather*}
and for
\begin{gather*}
{ \rm III_0} \quad w= -1+O\big(x^{-1}\big),\qquad x\to\infty \qquad \text{for} \quad \delta=2.
 \end{gather*}

\section[Truncated solutions in the family ${\rm I}_0$]{Truncated solutions in the family $\boldsymbol{{\rm I}_0}$}

In this section we state the main results; their proofs are found in Section~\ref{Pffam1}.

 \begin{Theorem}\label{TheoremFam1} Assume $\alpha\beta\gamma\ne 0$. Let $w(x)$ be a solution of ${\rm P}_{\rm V}\big(\alpha,\beta,\gamma,-\tfrac12\big)$ such that there exists some $\phi\in \big({-}\tfrac{\pi}2,\tfrac{\pi}2\big)$ so that
 \begin{gather}\label{alongphi}
 w(x)= \sqrt{-2\beta} x^{-1} (1+o(1))\qquad \text{as} \quad x\to\infty \quad \text{along} \quad \arg x=-\phi.
 \end{gather}

 Then \eqref{alongphi} holds along any $\phi\in \big({-}\tfrac{\pi}2,\tfrac{\pi}2\big)$.
\begin{enumerate}\itemsep=0pt
 \item[$(i)$] Furthermore, $w(x)$ is asymptotic to a unique power series solution:
 \begin{gather*}
 w(x)\sim \tilde{w}_0(x)= \sum_{n=1}^\infty {w}_{0;n}x^{-n}\qquad \text{as} \quad x\to\infty \quad \text{along any} \quad \arg x=-\phi,
 \end{gather*}
 where $|\phi|<\tfrac{\pi}2$, ${w}_{0;1}=\sqrt{-2\beta}$.

 \item[$(ii)$] The complete formal solution $($transseries$)$ along any half-line $\arg x=-\phi$ with $0<|\phi|<\tfrac{\pi}2$ is
 \begin{gather}\label{tranI0}
 \tilde{w}(x)=\tilde{w}_0(x)+\sum_{k=1}^{\infty}C^k{\rm e}^{-kx}x^{-qk}\tilde{w}_k(x)\qquad \text{where} \quad \tilde{w}_k(x)=\sum_{n= 0}^{\infty} w_{n;k}x^{-n},
 \end{gather}
 $C$ is an arbitrary constant and
 \begin{gather*}
 q=\gamma+2\sqrt{-2 \beta}.
 \end{gather*}

 \item[$(iii)$] $w(x)$ has a Borel summed transseries representation: there exist $C_\pm$, where $C_+-C_-$ is a~multiple of the Stokes constant, such that for any $\epsilon>0$ there is $R_\epsilon>0$
 \begin{gather}\label{sumtrans}
 w(x)= \begin{cases}\displaystyle w_0(x)+\sum_{k=1}^{\infty}C_+^k{\rm e}^{-kx}x^{-qk}w_{+;k}(x) & \text{for} \ \arg x\in \big(0,\frac{\pi}2-\epsilon\big),\ |x|>R_\epsilon,\\
 \displaystyle w_0(x)+\sum_{k=1}^{\infty}C_-^k{\rm e}^{-kx}x^{-qk}w_{-;k}(x) & \text{for} \ \arg x\in \big({-}\frac{\pi}2+\epsilon,0\big),\ |x|>R_\epsilon,
 \end{cases}\!\!\!\!
 \end{gather}
where $w_0(x)$, $w_{\pm;k}(x)$ are the Borel sums of $\tilde{w}_0(x)$, $\tilde{w}_k(x)$ along half-lines ${\rm e}^{{\rm i}\phi}\RR_+$ $($where $\phi=-\arg x)$. For $\arg x=0$ the transseries is \'Ecalle--Borel summable.\footnote{It is obtained by special averages of Laplace transforms in the upper and the lower half plane~\cite{Duke}.}
The function series~\eqref{sumtrans} converge for $x$ satisfying $|x|>R$ and $|C_\pm {\rm e}^{-x}x^{-q}|<\mu$ for suitable $\mu>0$ and $R>0$.
 In particular, if $\operatorname{Re} q>0$ then $w(x)$ is analytic in the right half plane for $|x|>R$ $($for some $R>0)$.

\item[$(iv)$] If $\operatorname{Re} q>0$ the constants $C_\pm$ can be determined from
 \begin{gather*}
 w(x)\sim \sum_{n=1}^{\lfloor\operatorname{Re} q\rfloor} {w}_{0;n} x^{-n}+C_+x^{-q} {\rm e}^{-x}\qquad \text{as} \quad x\to+{\rm i}\infty
 \end{gather*}
 and
 \begin{gather*}
 w(x)\sim \sum_{n=1}^{\lfloor \operatorname{Re} q\rfloor} {w}_{0;n} x^{-n}+C_-x^{-q} {\rm e}^{-x}\qquad \text{as} \quad x\to-{\rm i}\infty.
 \end{gather*}
 \end{enumerate}
\end{Theorem}

\begin{Corollary}[existence of tri-truncated solutions] Consider the unique truncated solu\-tion $w(x)$ as in Theorem~{\rm \ref{TheoremFam1}} with $C_+=0$. Then $w(x)$ is analytic for large $|x|$ in the left-half plane, for $\arg x\in\big({-}\tfrac{\pi}2+\epsilon, \tfrac{3\pi}2-\epsilon\big)$, $|x|>R_\epsilon$ $($for any $\epsilon>0)$.

Similarly, the unique truncated solution with $C_-=0$ is analytic for $|x|>R_\epsilon$ with $\arg x\in\big({-}\tfrac{3\pi}2+\epsilon,\tfrac{\pi}2-\epsilon\big)$. If $\operatorname{Re} q>0$ then we can take $\epsilon=0$.
\end{Corollary}

\begin{Remark} Truncated solutions in the left half-plane, satisfying
 \begin{gather*}
 w(x)= -\sqrt{-2\beta} x^{-1} (1+o(1))\qquad \big(x\to\infty \text{ along }\arg x=-\phi,\ |\phi+\pi|<\tfrac{\pi}2\big)
 \end{gather*}
form a one-parameter family, which, by Remark~\ref{Rem1} are obtained by replacing $(x,\gamma)$ with $(-x,-\gamma)$ in the representation given by Theorem~\ref{TheoremFam1}.
\end{Remark}

Truncated solutions $w(x)$ as in Theorem~\ref{TheoremFam1} with $C_+\ne 0$ develop arrays of poles. The following result shows the position of the array of poles closest to ${\rm i}\RR^+$.

\begin{Theorem}\label{arrays} Assume $\alpha\beta\gamma\ne 0$ and
 \begin{gather*}
 2 \alpha \ne \big(\sqrt {-2 \beta}-q-1\big)^2.
\end{gather*}

Let $w(x)$ be as in Theorem~{\rm \ref{TheoremFam1}} with nonzero constant $C_+$ in \eqref{sumtrans}. Them $w(x)$ has two arrays of poles located at
\begin{gather*}
x_{n;1,2}=2n\pi{\rm i}+\big(\gamma+2 \sqrt {-2 \beta}+2\big)\ln(2n\pi{\rm i})+\ln C_+-\ln \zeta_{1,2}+o(1), \qquad n\to+\infty,
\end{gather*}
where
 \begin{gather}\label{locsing}
 \zeta_{1,2}=\frac 2{\sqrt {-2 \beta}} \frac 1{\sqrt {-2 \beta}-q-1\pm \sqrt{2\alpha}}.
\end{gather}
\end{Theorem}

\section[Truncated solutions in the family ${\rm III}_0$]{Truncated solutions in the family $\boldsymbol{{\rm III}_0}$}

 In this section we state the main results; their proofs are found in Section~\ref{Pffam3}.

 \begin{Theorem}\label{TheoremFam3} Assume $\alpha\beta\gamma\ne 0$. Let $w(x)$ be a solution of ${\rm P}_{\rm V}(\alpha,\beta,\gamma,2)$ such that there exists some $\phi\in \big({-}\tfrac{\pi}2,\tfrac{\pi}2\big)$ so that
 \begin{gather}\label{asym1}
 w(x)=-1+o(1)\qquad \text{as} \quad x\to\infty \quad \text{along} \quad \arg x=-\phi.
 \end{gather}

 Then \eqref{asym1} holds along any $\phi\in \big({-}\tfrac{\pi}2,\tfrac{\pi}2\big)$.
\begin{enumerate}\itemsep=0pt
 \item[$(i)$] Furthermore, $w(x)$ is asymptotic to a unique power series solution:
 \begin{gather*}
 w(x)\sim -1+\tilde{w}_0(x)=-1+ \sum_{n=1}^\infty {w}_{0;n}x^{-n}\qquad x\to\infty \quad \text{along any} \quad \arg x=-\phi,
 \end{gather*}
where $|\phi|<\tfrac{\pi}2$.
 \item[$(ii)$] The complete formal series solution $($transseries$)$ along any half-line in the right half-plane has the form
 \begin{gather*}
 \tilde{w}(x)=\tilde{w}_0(x)+\sum_{k=1}^{\infty}C^k{\rm e}^{-kx}x^{-k/2}\tilde{w}_k(x),
 \end{gather*}
 where $C$ is an arbitrary constant, $\tilde{w}_k(x)$ are series in $x^{-n}$, $n\in\NN$.

 \item[$(iii)$] There are unique constants $C_\pm$ so that
 \begin{gather}\label{recoupCp}
 w(x)\sim -1+C_+ x^{-1/2} {\rm e}^{-x}\qquad \text{as} \quad x\to+{\rm i}\infty
 \end{gather}
 and
 \begin{gather*}
 w(x)\sim-1+C_- x^{-1/2} {\rm e}^{-x}\qquad \text{as} \quad x\to-{\rm i}\infty
 \end{gather*}
and $C_+-C_-$ is a constant which does not depend on the particular solution $($and it is a~multiple of the Stokes constant of the equation$)$.

 \item[$(iv)$] $w(x)$ is analytic in the right half plane for $|x|>R$ $($for some $R>0)$ and has the Borel summed transseries representation
 \begin{gather}\label{sumtrans2}
 w(x)= \begin{cases}\displaystyle w_0(x)+\sum_{k=1}^{\infty}C_+^k{\rm e}^{-kx}x^{-k/2}w_{+;k}(x) & \text{for} \ \arg x\in \big(0,\frac{\pi}2\big], \\
\displaystyle w_0(x)+\sum_{k=1}^{\infty}C_-^k{\rm e}^{-kx}x^{-k/2}w_{-;k}(x) & \text{for} \ \arg x\in \big[{-}\frac{\pi}2,0\big),
 \end{cases}
 \end{gather}
where $w_0(x)$, $w_{\pm;k}(x)$ are Laplace transforms on the half-lines ${\rm e}^{{\rm i}\phi}\RR_+$ $($where $\phi=-\arg x$ is between~$0$ and $\pm\tfrac{\pi}2)$ of the Borel transforms of the series $\tilde{w}_0(x)$, $\tilde{w}_k(x)$.
The domain of convergence of the series in~\eqref{sumtrans2} is the set of all $x$ with $|x|>R$ and $\big|C_\pm {\rm e}^{-x}x^{-1/2}\big|<\mu$ for some suitable $R,\mu>0$.
 \end{enumerate}
\end{Theorem}

\begin{Corollary}[existence of tri-truncated solutions] Consider the unique truncated solution $w(x)$ as in Theorem~{\rm \ref{TheoremFam3}} with $C_+=0$. Then $w(x)$ is analytic for large $|x|$ also in the left-half plane, for $\arg x\in\big({-}\tfrac{\pi}2+\epsilon, \tfrac{3\pi}2-\epsilon\big)$.

Similarly, the unique truncated solution with $C_-=0$ is analytic for large $|x|$ with $\arg x\in\big( {-}\tfrac{3\pi}2+\epsilon,\tfrac{\pi}2-\epsilon\big)$.
\end{Corollary}

\begin{Remark} Truncated solutions in the left half-plane, satisfying
 \begin{gather*}
 w(x)= -1+o(1)\qquad \big(x\to\infty\text{ along }\arg x=-\phi\text{ with }0<|\phi+\pi|<\tfrac{\pi}2 \big)
 \end{gather*}
form a one-parameter family, which, by Remark~\ref{Rem1} are obtained by replacing $(x,\gamma)$ with $(-x,-\gamma)$ in the representation given by Theorem~\ref{TheoremFam3}.
\end{Remark}

Truncated solutions $w(x)$ as in Theorem~\ref{TheoremFam3} with $C_+\ne 0$ have arrays of poles near ${\rm i}\RR_+$. The following result shows the position of the closest array of poles.

\begin{Theorem}\label{arrays2} Assume $\alpha\beta\gamma\ne 0$. Let $w(x)$ be as in Theorem~{\rm \ref{TheoremFam3}} with nonzero constant $C_+$ in~\eqref{recoupCp} and~\eqref{sumtrans2}. Them $w(x)$ has an array of poles located at
\begin{gather*}
x_n=2n\pi {\rm i} -\frac12\ln(2n\pi {\rm i}) -\frac{{\rm i}\pi}{2}-\ln(-C_+/2)+o(1),\qquad n\to\infty.
\end{gather*}
\end{Theorem}

We {\em conjecture} that there are, in fact, two poles near each~$x_n$.

\section[Proofs for family ${\rm I}_0$]{Proofs for family $\boldsymbol{{\rm I}_0}$}\label{Pffam1}

\subsection{Proof of Theorem \ref{TheoremFam1}}\label{PfTh2}

Denote
\begin{gather*}
m=\sqrt{-2 \beta},\qquad q=\gamma+2m.
\end{gather*}
Then $w(x)$ satisfies
\begin{gather}\label{eqP5}
w'' = {\frac {( 3 w-1) w'^{2}}{ 2 w (w-1) }} -{\frac {w' }{x}} +{\frac {\alpha w (w-1) ^{2}}{{x}^{2}}}-\frac{{m}^{2}}2 {\frac {( w -1) ^{2}}{{x}^{2}w }}+{\frac {(q -2 m) w }{x}} -\frac12 {\frac{w (w +1) }{w -1}}
\end{gather}
for which we study solutions satisfying $w=mx^{-1}+O\big(x^{-2}\big)$ as $x\to \infty$.

To normalize the equation, substitute
\begin{gather}\label{subyu}
w(x)=\frac mx\left(1-\frac qx +u(x)\right),
\end{gather}
which transforms \eqref{eqP5} into
\begin{gather}\label{normP5}
u''=\left(1+\frac{2q}x\right)u+f\big(x^{-1},u,u'\big),
\end{gather}
where $f$ is analytic at $(0,0,0)$ and $f=O\big(x^{-2}\big)+O\big(u^2\big)+O\big(u'^2\big)+O(uu')$. To show that this a normal form for second order equations, we turn it into a first order system by substituting
\begin{gather}\label{subuy}
 u (x) =y_1(x) +y_2(x) ,\qquad u' (x) = \left( -1-{\frac {q}{x }} \right) y_1(x) + \left( 1+{\frac {q}{x}} \right) y_2(x),
 \end{gather}
upon which \eqref{normP5} is turned into a first order system in normal form \eqref{normform}:
\begin{gather}\label{sysy}
\frac{{\rm d}}{{\rm d}x}\left[ \begin{matrix} y_1\\ y_2\end{matrix}\right]=
 \left[ \begin{matrix}
 -1-{\dfrac {q}{x}} & 0
 \\ 0 & 1+{\dfrac {q}{x}}
\end{matrix}\right] \left[ \begin{matrix} y_1\\ y_2\end{matrix}\right]
 +g\big(x^{-1},y_1,y_2\big)\left[ \begin{matrix} 1\\ -1
 \end{matrix}\right],
 \end{gather}
 where
 \begin{gather}
 g\big(x^{-1},y_1,y_2\big)=\frac 12\frac x{x+q} \left[- f\left(x^{-1},y_1+y_2, \left(1+\frac qx\right)(y_2-y_1)\right) \right.\nonumber\\
\left. \hphantom{g\big(x^{-1},y_1,y_2\big)=}{} + \frac{q(q+1)}{x^2} y_1+ \frac{q(q-1)}{x^2} y_2 \right].\label{formg}
 \end{gather}

 There are general theorems that can be applied for differential equations in normal form; these are presented, for convenience, in the Appendix~\ref{AppB}: applying Theorem~\ref{SumTransS} to the system \eqref{sysy}, \eqref{formg} and then reverting the substitutions \eqref{subuy}, \eqref{subyu}, Theorem~\ref{TheoremFam1} follows.

 \subsection{Proof of Theorem \ref{arrays}}

Searching for asymptotic expansions of the form \eqref{asser} setting $\xi=C{\rm e}^{-x}x^{-q}$, plugging in an asymptotic series $u(x)\sim F_0(\xi)+\tfrac 1x F_1(\xi)+\tfrac 1{x^2} F_2(\xi)+\cdots$ in \eqref{normP5} and expanding under the assumption that $\xi\gg x^{-k}$ for all $k$ we obtain that all $F_n$ are polynomials:
 \begin{gather*}
 F_0(\xi)=\xi,\qquad F_1(\xi)=c_1\xi,\\
 F_2(\xi)=m ( m-q-1 ) {\xi}^{2}+c_2 \xi-\tfrac 12 {m}^{2}+\tfrac 32 {q}^{2}-a+\tfrac12 ,\qquad \ldots,
\end{gather*}
where $c_1$, $c_2$ are uniquely determined in terms of the parameters. By Theorem~\ref{Thinvent} we have
\begin{gather}\label{asseru}
u(x)\sim \sum_{m=0}^{\infty} x^{-m}{F}_m(\xi(x))\qquad \text{for} \ x\to\infty \ \text{with} \
\arg x\in \big[{-}\tfrac{\pi}2+\delta, \tfrac{\pi}2-\delta\big],\ |\xi(x)|<\delta_1
\end{gather}
for some $\delta_1>0$.

On the other hand, $F_0$ has no singularities, so Theorem~\ref{Tsing} yields no additional information on the position of the first array of poles, beyond the right half plane.

 A direct calculation of the first few terms of the transseries \eqref{tranI0} unveils a non-generic structure, namely that $w_{n,k}=0$ for all $k=0,1,\ldots,n-1$ (at least for $n=1,\dots,4$). This peculiar structure suggests (by formal re-arrangement of the transseries) to use instead the second scale
 \begin{gather*}
 \zeta=C{\rm e}^{-x}x^{-q-2}
 \end{gather*}
 and look for an expansion of the form
 \begin{gather}\label{formFser}
 u(x)\sim x^2\Phi (\zeta)+\sum_{n=-1}^\infty \frac 1{x^n}\Phi_n (\zeta), \qquad \zeta(x)\gg x^{-k}\ \text{for all} \ k,
 \end{gather}
where $\Phi (\zeta)=\zeta+O\big(\zeta^2\big)$ when $\zeta\to 0$ and $\Phi_n $ are analytic at $\zeta=0$, and for $x$ so that $|\zeta|=\big|C{\rm e}^{-x}x^{-q-2}\big|<\mu$, where $\mu$ is small enough (to be determined).

 Introducing the formal expansion \eqref{formFser} in \eqref{normP5} and expanding in powers of $x^{-1}$, we obtain that $\Phi (\zeta)$ must satisfy
 \begin{gather*}
 {\zeta}^{2} \Phi ''+ \zeta \Phi ' +\frac12 \Phi = \alpha {m}^{2} \Phi ^{3} +\frac 32\frac{\zeta^2 \Phi '^2}{\Phi }
 \end{gather*}
having the general solution
\begin{gather}\label{genF}
\Phi (\zeta)=\frac{2\zeta/C_1}{8\alpha m^2/C_1^2-(\zeta-C_2)^2},
\end{gather}
where $C_{1,2}$ are determined from the condition that $\Phi (\zeta)=\zeta+O\big(\zeta^2\big)$ and that $\Phi_{-1}(\zeta)$ be analytic at $\zeta=0$, yielding
 \begin{gather}\label{valC1C2}
 C_1=-2 \beta \big[2 \alpha -(m-q-1)^2\big],\qquad C_2=-\frac 2m \frac {m-q-1}{ 2 \alpha -(m-q-1)^2 }.
\end{gather}

With the notation \eqref{locsing} (recall that $m=\sqrt{-2\beta}$) and using \eqref{genF}, then \eqref{valC1C2} becomes
 \begin{gather*}
 \Phi(\zeta)=\frac {\zeta \zeta_1\zeta_2}{(\zeta-\zeta_1) (\zeta-\zeta_2)}.
 \end{gather*}
Letting $\mu<\min |\zeta_{1,2}|$, the next term in the expansion is
\begin{gather}
\Phi_{-1}(\zeta)={\frac { \zeta ( 2 \zeta-\zeta_{{1}}-\zeta_{{2}} ) C_{{3}}}{ ( \zeta-\zeta_{{2}} ) ^{2} ( \zeta-\zeta_{{1}} ) ^{2}}}+{\frac {
 \zeta \big( {\zeta}^{2}-\zeta_{{1}}\zeta_{{2}} \big) C_{{4}}}{ ( \zeta-\zeta_{{2}} ) ^{2} ( \zeta-\zeta_{{1}} ) ^{2}}}\label{Phim1}\\
\hphantom{\Phi_{-1}(\zeta)=}{} -\frac1{2a} {\frac { \zeta \big( \zeta_{{1}}\zeta_{{2}}\sqrt {2} ( \zeta_{{1}}+\zeta_{{2}} ) {a}^{3/2}+ ( \zeta_{{1}}-\zeta_{{2}} ) \big( {-}\frac 1{\sqrt {2} }\left( \zeta_{{1}}-\zeta_{{2}} \right) \sqrt {a}+a\zeta_{{2}}\zeta_{{1}} \big) \big) ( \zeta_{{1}}-\zeta_{{2}} ) }{ ( \zeta-\zeta_{{2}} ) ^{2} ( \zeta-\zeta_{{1}} ) ^{2}}},\nonumber
\end{gather}
where $C_{3,4}$ are determined from the condition that the next term, $\Phi_{0}$, be analytic at $\zeta=0$.

To justify that the actual truncated solution $w(x)$ is also singular near $\zeta_1$ and near $\zeta_2$, we first note that, since $\zeta=\xi x^{-2}$ and $u(x)\sim x^2\Phi$ then
\begin{gather*}
x^2\Phi(\zeta)=\frac{\xi}{\big( 1-x^{-2}\xi\zeta_1^{-1}\big)\big( 1-x^{-2}\xi\zeta_2^{-1}\big)}=\xi+\frac 1{x^2}\xi\big(\zeta_1^{-1}+\zeta_2^{-1}\big)+\cdots,
\end{gather*}
which is a convergent series in powers of $x^{-2}$. Also $x\Phi_{-1}(\zeta)$ has a similar convergent expansion.

Therefore $x^2\Phi(\zeta)=F_0(\xi)+O\big(x^{-2}\big)$ and, in view of \eqref{Phim1}, $x\Phi_{-1}(\zeta)=O\big(x^{-1}\big)$ therefore \eqref{asseru} implies that
\begin{gather*}
u(x)\sim x^2\Phi(\zeta)+O\big(x^{-1}\big)
\end{gather*}
for $x\to\infty$ in the same region where \eqref{asseru} holds (and $|x|>1$, $|\zeta|<\mu$). The same argument as in \cite[Section~4.6]{Invent} implies that $u(x)$ has singularities within~$o(1)$ distance of the singularities of~$\Phi(\zeta)$.

\section[Proofs for the family ${\rm III}_0$]{Proofs for the family $\boldsymbol{{\rm III}_0}$}\label{Pffam3}

The existence of a unique power series formal solution $\tilde{w_0}(x)=-1+o(1)$ is established by standard techniques.

It is also relatively algorithmic to obtain the form of the transseries. Since the procedure may not be well known some details are provided here, also illustrating why $\delta=2$ is a natural choice.

Plugging in a formal solution of the type $\tilde{w}(x)=\tilde{w_0}(x)+\epsilon g(x)+O\big(\epsilon^2\big)$ in the equation \eqref{P5} and expanding in power series in $\epsilon$, the coefficient of $\epsilon^0$ vanishes, since $\tilde{w_0}(x)$ is already a formal solution. Next, the coefficient of $\epsilon$ is a linear second order differential equation for $g(x)$, with coefficients given in terms of $\tilde{w_0}(x)$ and its derivatives. It has the form
\begin{gather*}
g''(x)+\frac 1x g'(x)-\frac {\delta}2 g(x)=\frac 1{x^2} R\big(x^{-1},g,g'\big)
\end{gather*}
with two independent solutions
\begin{gather*}
g_{\pm}(x)=x^{-1/2}\exp\big(\pm x\sqrt{\delta/2}\big) \big(1+O\big(x^{-1}\big)\big),\end{gather*}
where we see that it is convenient to take $\delta=2$.

Let $w(x)$ satisfy the assumptions of Theorem~\ref{TheoremFam3}. To normalize the equation, substitute
\begin{gather*}
w(x)=-1+\frac {\gamma}x+u(x),
\end{gather*}
so that $u(x)=O\big(x^{-2}\big)$ and satisfies an equation of the form
\begin{gather*}
u''=\left(1-\frac{1}x\right)u+f\big(x^{-1},u,u'\big),
\end{gather*}
where $f$ is analytic at $(0,0,0)$ and $f=O\big(x^{-2}\big)+O\big(u^2\big)+O\big(u'^2\big)+O(uu')$. This is a normal form for a second order equation, by the argument in Section~\ref{PfTh2}: the normalizing transformation is \eqref{subuy} with $q=-1/2$ and the normal form as a first order system is of the type \eqref{sysy} with $q=-1/2$. Theorem~\ref{SumTransS} applies, and reverting the substitutions, Theorem~\ref{TheoremFam3} follows.

\subsection{The first array of poles: proof of Theorem~\ref{arrays2}}
Setting $\xi= {C{{\rm e}^{-x}}} {x}^{-1/2}$, plugging in $\mathcal{P}(\alpha,\beta,\gamma,2)$ an expansion~\eqref{asser} and expanding, it follows that $F_0(\xi)$ must satisfy
\begin{gather*}
{\xi}^{2} F_0'' +\xi F_0' = {\frac {{\xi}^{2} ( 3 F_0 -4 ) F _0'^{2}}{2 ( F_0 -2 ) ( F_0 -1 ) }} +2 {\frac {F_0 ( F_0 -1 ) }{F_0 -2}},
\end{gather*}
whose unique solution analytic at $\xi=0$ with $F_0(\xi)=\xi+O\big(\xi^2\big)$, $\xi\to 0$, is
\begin{gather*}
F_0  ( \xi ) ={\frac {\xi}{( \xi/4+1) ^{2}}}.
\end{gather*}

Let $x_n$ be solutions of $\xi(x)/4+1=0$. By Theorem~\ref{Tsing}, if $x_n$ solve $C{\rm e}^{-x}x^{-1/2}=-4$, then solutions $u(x)$ have singularities located at $x_n+o(1)$ for large $n$ and Theorem~\ref{arrays2} follows.

\begin{Remark} While $F_0(\xi)$ has a pole of order two, it is known that the poles of $P_V$ have order one if $\alpha\ne 0$. This suggests that there are two array of poles within $O\big(n^{-1}\big)$ distance from each other, as the expansion in~$x^{-1}$ would collapse them into a double pole. The author has no rigorous argument at this time that this is indeed the case, and it is formulated here as a~conjecture.
\end{Remark}

\appendix

\section{Results used}\label{AppB}

\subsection{Summation of transseries formal solutions}

\begin{Theorem}[\cite{Duke}]\label{SumTransD} Consider the (nonlinear) system of first order differential equations:
\begin{gather}\label{normform}
\bfy'+\left(\Lambda-\frac{1}{x}A\right)\bfy=\mathbf{g}(x^{-1},\bfy),\qquad \bfy\in\RR^d,
\end{gather} with $\mathbf{g}=O\big(x^{-2}\big)+O\big(|\bfy|^2\big)$ for $x\to\infty$, $|\bfy|\to 0$
and $\Lambda={\rm diag}(\lambda_1,\ldots,\lambda_d)$, $A={\rm diag}(\alpha_1,\ldots,\alpha_d)$.
Assume the following non-resonance condition: any collection among $\lambda_1,\ldots,\lambda_d$ which are in the same open half-plane are linearly independent over~$\ZZ$.

Assume, for simplicity, that the independent variable, $x$, is scaled so that $\lambda_1=1$.
\begin{enumerate}\itemsep=0pt
\item[$(i)$] Then \eqref{normform} has formal solutions
\begin{gather}\label{gtrans}
\tilde{\bfy}=\tilde{\bfy}(x;\mathbf{C})= \tilde{\bfy}_\mathbf{0}(x)+\sum_ {\mathbf{k}\in\NN^d\setminus\mathbf{0}} \mathbf{C}^\mathbf{k}{\rm{e}}^{-\boldsymbol\lambda \cdot \mathbf{k}x} \tilde{\bfy}_\mathbf{k}(x) \qquad \text{with} \quad \tilde{\bfy}_\mathbf{k}(x)=x^{\boldsymbol\alpha\cdot \mathbf{k}}\tilde{\bfs}_\mathbf{k}(x),
\end{gather}
where $\tilde{\bfs}_\mathbf{k}(x)$ are integer power series in $x^{-1}$.
A formal solution \eqref{gtrans} is a transseries for $x\to\infty$ along any direction along which all the exponentials present are decaying, i.e., along any direction in the sector
\begin{gather}\label{Strans}
S_{\rm trans}=\{x\in\CC \,|\, \operatorname{Re}(\lambda_j x)>0\ \text{for all} \ j\ \text{with} \ C_j\ne 0\} .
\end{gather}

\item[$(ii)$] A transseries is Borel summable for large $x$ along any direction of argument $\phi$, where the sector $a_1<\arg x<a_2$ contains only one Stokes line, $\arg x=0$. The Borel sum is an actual solution: for any $\epsilon>0$ there is $R_\epsilon>0$ so that
 \begin{gather}\label{GenTrans}
 \bfy(x)=
 \begin{cases}
\displaystyle\mathcal{L}_\phi \mathbf{Y}_\mathbf{0} (x)+ \!\sum_ {\mathbf{k}\in\NN^d\setminus\mathbf{0}}\! \mathbf{C_+}^\mathbf{k}{\rm{e}}^{-\boldsymbol\lambda \cdot \mathbf{k}x} \mathcal{L}_\phi\mathbf{Y}_\mathbf{k}(x),\
 - \phi= \arg x\in(0,a_2\!-\!\epsilon),\ |x|>R_\epsilon,\\
 \displaystyle \mathcal{L}_\phi \mathbf{Y}_\mathbf{0} (x)+ \!\sum_ {\mathbf{k}\in\NN^d\setminus\mathbf{0}}\! \mathbf{C_-}^\mathbf{k}{\rm{e}}^{-\boldsymbol\lambda \cdot \mathbf{k}x} \mathcal{L}_\phi\mathbf{Y}_\mathbf{k}(x),\
 - \phi= \arg x\in(a_1\!+\!\epsilon,0),\ |x|>R_\epsilon,
 \end{cases}\hspace{-20mm}
 \end{gather}
where $\mathbf{Y}_\mathbf{k}=\mathcal{B}_\phi\tilde{\bfy}_\mathbf{k}$ $($the analytic continuation of the Borel transform of $\tilde{\bfy}_\mathbf{k}$ along the direction of argument $\phi)$.
Along $\arg x=0$ balanced averages of $\mathcal{L}_\phi\mathbf{Y}_\mathbf{k}$ sum to the solution $\bfy(x)$ $($\'Ecalle--Borel summation$)$.
The solution $\bfy(x)$ is analytic for $|x|>R_\epsilon$ with $\arg x\in(a_1+\epsilon,a_2-\epsilon)$.
Conversely, any solution asymptotic to $\tilde{\bfy}_0$ for $x\to+\infty$ has a~representation \eqref{GenTrans}.

\item[$(iii)$] Only the first component $C_1$ of the constant beyond all orders in~\eqref{GenTrans} changes when $\arg x$ crosses the Stokes line $\arg x=0$, corresponding to $\lambda_1=1$. The change of this constant depends only on the equation: it is a multiple of the first Stokes constant.
\end{enumerate}
\end{Theorem}

 \begin{Note} More recent results showed that the region of convergence of \eqref{GenTrans} is in fact given by conditions of the form $|x|>R$ and $\big|C_i{\rm e}^{-\lambda_ix}x^{\alpha_i}\big|<\mu_i$, $i=1,\ldots,d$ for suitable constants $\mu_i$ and $R$ ($C_i=0$ if $\operatorname{Re}\lambda_i> 0$) \cite{Invent}, see also \cite{OBook}.
 \end{Note}

We use Theorem~\ref{SumTransD} in the particular case when $d=2$, $\lambda_1=1$, $\lambda_2=-1$, in which case the sectors~\eqref{Strans} become the right (respectively left) half plane and the constants have the form $ \mathbf{C}=(C,0)$ (respectively $ \mathbf{C}=(0,C)$); we also assume that solutions have only poles as moving singularities, as it is the case for the Painlev\'e equations (otherwise domains of analyticity would have to extend on Riemann sheets). Theorem~\ref{SumTransD} takes the following simpler formulation:

\begin{Theorem}\label{SumTransS}Consider the two-dimensional system of first order differential equations~\eqref{normform} with $\Lambda={\rm diag}(1,-1)$, $A={\rm diag}(\alpha_1,\alpha_2)$. Then:
\begin{enumerate}\itemsep=0pt
\item[$(i)$] The system has a one parameter family of transseries solutions
\begin{gather}\label{gtrans1}
\tilde{\bfy}=\tilde{\bfy}(x;C)= \tilde{\bfy}_{0}(x)+\sum_ {k\geq 1} {C}^{k}{\rm{e}}^{-{k}x} \tilde{\mathbf{y}}_{k}(x),\ C\in\CC, \ \text{along} \ \arg x=-\phi, \ |\phi|<\tfrac {\pi}2,
\end{gather}
where $\tilde{\bfy}_\mathbf{k}(x) =x^{{k}\alpha_1} \tilde{\mathbf{s}}_\mathbf{k}(x)=x^{{k}\alpha_1} \sum\limits_{n=0}^\infty\tilde{\mathbf{s}}_\mathbf{k,n}x^{-n}\text{ and }\tilde{\bfy}_{0}(x)=O\big(x^{-2}\big)$.
\item[$(ii)$] For any $C$ the formal solution \eqref{gtrans1} is Borel summable for large $x$ along any direction of argument $\phi=-\arg x$ if $0<|\phi|<\tfrac{\pi}2$. $($Along $\arg x=0$ it is \'Ecalle--Borel summable.$)$ The Borel sum is an actual solution: for any constant $C_+$, and any $\epsilon>0$ there exists $C_-$ $(C_+-C_-$ is a multiple of the Stokes constant, {\rm \cite{P1us})} and $R_\epsilon>0$ such that
 \begin{gather}\label{summedtrans}
 \bfy(x)=
 \begin{cases}
\displaystyle \mathcal{L}_{\phi}\mathbf{Y}_\mathbf{0} (x)+ \sum_ {k\geq 1} C_+^{k}{\rm{e}}^{-{k}x} \mathcal{L}_{\phi}\mathbf{Y}_{k}(x), & - \phi= \arg x\in\big(0,\frac{\pi}2-\epsilon\big),\ |x|>R_\epsilon, \\
\displaystyle \mathcal{L}_{\phi}\mathbf{Y}_\mathbf{0} (x)+ \sum_ {k\geq 1} C_-^{k}{\rm{e}}^{-{k}x} \mathcal{L}_{\phi}\mathbf{Y}_{k}(x),& - \phi= \arg x\in\big({-}\frac{\pi}2+\epsilon,0\big),\ |x|>R_\epsilon.
 \end{cases}\hspace{-10mm}
 \end{gather}
The domain of convergence of the function series \eqref{summedtrans} given by
\begin{gather*}
\big\{x \,|\, |x|>R,\, \big|C_\pm {\rm e}^{-x}x^{\alpha_1}\big|<\mu\big\}
\end{gather*} for $\mu>0$ small enough and $R>0$ large enough.
Conversely, any solution asymptotic to~$\tilde{\bfy}_0$ for $x\to+\infty$ has a representation \eqref{summedtrans} for some $C_\pm$.

\item[$(ii)$] Similarly, the system has a transseries solution
\begin{gather}\label{gtrans2}
\tilde{\bfy}=\tilde{\bfy}(x;C)= \tilde{\bfy}_{0}(x)+\sum_ {k\geq 1} {C}^{k}{\rm{e}}^{{k}x} \tilde{\bfy}_{k}(x) \qquad \text{for} \ \arg x=-\phi,\ 0<|\phi+\pi|<\frac {\pi}2,
\end{gather}
 where $ \tilde{\bfy}_\mathbf{k}(x)=x^{{k}\alpha_2}\tilde{\bfs}_\mathbf{k}(x)$ with $\tilde{\bfs}_\mathbf{k}(x)$ are integer power series.
\end{enumerate}

 Statements similar to those in $(i)$ hold regarding Borel summability \eqref{gtrans2} along any direction of argument $\arg x\ne\pi$ in the left half plane, and \'Ecalle--Borel summability along $\arg x=\pi$, the sum being an actual solution, analytic in a domain $|x|>R$ and $\big|C_\pm {\rm e}^{x}x^{\alpha_2}\big|<\mu$ for some $\mu>0$ and $R>0$.
 Conversely, any solutions asymptotic to $\tilde{\bfy}_0$ for $x\to-\infty$ is the Borel sum of such a~transseries.
 \end{Theorem}

 \subsection{Arrays of singularities bordering the sector of analyticity}

Theorem~\ref{SumTransD} establishes existence of solutions, in one-to-one correspondence with formal trans\=se\-ries solutions, and which are analytic for large~$x$ in the sector~\eqref{Strans} where these formal solutions are defined (i.e., they are well-ordered with respect to~$\gg$). In~\cite{Invent} it is further shown that on the boundary of the sector of analyticity, these solutions develop arrays of singularities.

In the particular case when $d=2$, $\lambda_1=1$, $\lambda_2=-1$ of this paper the results in \cite{Invent} are as follows.
Denote
\begin{gather*}
\xi=\xi(x)=C{\rm e}^{-x}x^{\alpha_1}.
\end{gather*}
For $x$ near ${\rm i}\RR_+$, i.e., when $\xi\gg x^{-k}$ for all $k>0$, the transseries \eqref{gtrans1} can be formally reordered as
\begin{gather*}
\sum_ {k\geq 0} \xi^k \tilde{\mathbf{s}}_{k,0}(x) + \frac 1x\sum_ {k\geq 0} \xi^k \tilde{\mathbf{s}}_{k,1}+\frac 1{x^2}\sum_ {k\geq 0} \xi^k \tilde{\mathbf{s}}_{k,2}+\cdots.
\end{gather*}

It turns out that the series in $\xi$ are convergent, and the resulting expansion is asymptotic to the solution the original transseries summed to:

\begin{Theorem}[\cite{Invent}]\label{Thinvent} Let $\delta, c>0$. There exists $\delta_1>0$ so that for $|\xi|<\delta_1$ the power series
\begin{gather*}
\mathbf{F}_m(\xi)=\sum_ {k=0} ^{\infty} \xi^k \tilde{\mathbf{s}}_{k,m},\qquad m=0,1,2,\ldots
\end{gather*}
converge.
Furthermore
\begin{gather}\label{asser}\mathbf{y}(x)\sim \sum_{m=0}^{\infty} x^{-m}\mathbf{F}_m(\xi(x))\qquad \text{for} \ x\to\infty \ \text{with} \
\arg x\in \big[{-}\tfrac{\pi}2+\delta, \tfrac{\pi}2-\delta\big],\ |\xi(x)|<\delta_1.\!\!\!\!\!
\end{gather} The asymptotic series is uniform, it is differentiable and satisfies Gevery-like estimates.
\end{Theorem}

\begin{Note} $\mathbf{F}_0(0)=0$ and $\mathbf{F}'_0(0)=1$.\end{Note}

In fact \eqref{asser} is valid in a larger domain, up to distance $o(1)$ of the singularities of $\mathbf{F}_0$ ($\mathbf{F}_m$~with $m>0$ can have no other singularities). In general $\mathbf{F}_0$ has branch point singularities, and a Riemann surface needs to be considered. But Painlev\'e equations have no movable branch points, so the singularities of~$\mathbf{F}_0$ can only be poles. Then for simplicity we state here the general result of~\cite{Invent} in this case only.

Let $\rho_{1,2}$ so that the small term $\mathbf{g}$ in \eqref{normform} is analytic in the polydisk $\big|x^{-1}\big|<\rho_1$, $|\mathbf{y}|<\rho_2$.

Let $\Xi$ be a finite set (possibly empty) of poles of $\mathbf{F}_0$. Let $\mathcal{D}\subset\CC\setminus \Xi$ be open, connected, relatively compact, containing $|\xi|<\delta_1$, so that $\mathbf{F}_0$ is analytic in an $\epsilon$-neighborhood of $\mathcal{D}$ ($\epsilon>0$), so that $\sup\limits_\mathcal{D} |\mathbf{F}_0(\xi)|=\rho_3<\rho_2 $.

We need the counterpart of $\mathcal{D}$ in the $x$-plane. Let
\begin{gather*} X=\xi^{-1}(\mathcal{D})\cap \big\{ |x|>R, \,\arg x\in \big[{-}\tfrac{\pi}2+\delta, \tfrac{\pi}2-\delta\big]\big\}.\end{gather*}

\begin{Theorem}[\cite{Invent}]\label{Tsing}All $\mathbf{F}_m$ with $m\geq 1$ are analytic on $\mathcal{D}$ and for $R$ large enough the asymptotic expansion~\eqref{asser} holds for $x\to\infty$ with $x\in X$.

Furthermore, assume that $\mathbf{F}_0$ is singular at $\xi_s\in \mathcal{D}$. Then $\mathbf{y}(x)$ is singular near points $x$ with $\xi(x)=\xi_s$, more precisely at
\begin{gather*}
x_n=2n\pi {\rm i}+\alpha_1\ln(2n\pi {\rm i})+\ln C-\ln\xi_s+o(1)\qquad \text{as} \quad n\to\infty.
\end{gather*}
\end{Theorem}

\subsection*{Acknowledgements}

The author is grateful to the editors for valuable references and information, and to the referees careful reading of the manuscript and for their helpful comments.

\pdfbookmark[1]{References}{ref}
\LastPageEnding

\end{document}